\documentclass[a4paper,10pt]{article}

\def \Z {{\mathbf {Z}}}
\def \R {{\mathbf {R}}}

\def\uu{\bigsqcup}
\def\eps{\varepsilon}
\title{Эргодические  гомоклинические группы, сидоновские  конструкции и пуассоновские надстройки }
\author{В.В. Рыжиков}
\date{}

\usepackage[T1]{fontenc} 
\usepackage[cp1251]{inputenc} 
\usepackage[russian]{babel}
\begin{document}
\Large
\maketitle

\begin{abstract}{\Large Предложены новые  примеры перемешивающих  преобразований на пространстве с бесконечной мерой: так называемые сидоновские конструкции ранга 1.  Для класса бесконечных преобразований  получено быстрое убывание корреляций, недавно обнаруженное А.А. Приходько для динамических систем с простым спектром, действующих в на вероятностном пространстве.
Получен положительный ответ  на вопрос М.И. Гордина о существовании преобразований с нулевой энтропией и эргодическим гомоклиническим потоком. 
Рассмотрены модификации сидоновских конструкций, индуцирующие  пуассоновские надстройки с простым сингулярным спектром и гомоклиническим  бернуллиевским потоком. Дано новое доказательство теоремы Э. Руа о  кратном перемешивании   пуассоновских надстроек.
 
УДК  517.987.    MSC: Primary 28D05; Secondary 58F11.  

Ключевые слова:  Эргодические действия, конструкция ранга один,  множество Сидона, кратное перемешивание, пуассоновская надстройка,
гомоклинические преобразования, сингулярный спектр. 	}
\end{abstract}
%\begin{abstract}{\Large 

\section{Введение}
Эргодическая теория сохраняющих меру преобразований 
изучает  действия на вероятностном пространстве и действия на пространстве 
 с   сигма-конечной инвариантной мерой. Последние для  краткости будем называть бесконечными а    фазовое пространство  $(X,\mu)$, $\mu(X)=\infty,$ предполагается  изоморфным прямой с мерой Лебега.

Пуассоновская мера $\mu_\ast$
на пространстве конфигураций $X_\ast$ (см.\cite{I}, \cite{KSF},\cite{Ro})   
индуцирует непрерывное вложение группы $Aut(\mu)$  всех сохраняющих меру $\mu$ обратимых преобразований  в группу  $Aut(\mu_\ast)$.  Свойство перемешивания для  бесконечного преобразования $T$   означает
$$
\mu(T^nA\cap B) \to 0, \ n\to\infty,$$
для любых  множеств $A,B$  конечной меры.
Оно влечет за собой свойство перемешивания пуассоновской надстройки $T_\ast$:
$$
\mu_\ast(T_\ast^nV\cap W) \to \mu_\ast(V)\mu_\ast(W),  \ \ n\to\infty,$$
для  измеримых множеств $V,W\subset X_\ast$ $n\to\infty$.

Как мы  увидим, среди пуассоновских надстроек    только перемешивающие  обладают так называемыми 
гомоклиническими эргодическими группами. 
Следуя М.И. Гордину   \cite{G},  определим гомоклиническую группу
$H(T)$ преобразования   $T$: 
 $$H(T)=\{ S: \ T^nST^{-n}\to Id,  n\to\infty\}.$$

Если фазовое пространство $X$  имеет бесконечную меру, а преобразование $T:X\to X$ перемешивает,  то группа $F$ всех сохраняющих меру преобразований $S$ с носителем конечной меры, $\mu(supp(S))<\infty$,  лежит в группе $H(T)$.  
Группа $F_\ast$ лежит в $H(T_\ast)$ и является эргодической.

В \cite{G} установлен замечательный факт: эргодичность гомоклинической группы  автоморфизма  вероятностного пространства влечет за собой свойство перемешивания. Мы обобщим его  результат на  свойство  кратного перемешивания.  
 Это позволяет дать новое  доказательство теоремы Э.Руа о кратном перемешивании перемешивающих пуассоновских надстроек. 
Тем самым  в классе пуассоновских надстроек будет установлена 
эквивалентность свойств 
 перемешивания,
 кратного перемешивания всех порядков,
 эргодичности гомоклинической группы. 
Напомним, что перемешивание кратности 2 означает при $m,n\to\infty$
выполнение
$$
\mu_\ast(U\cap T_\ast^mV\cap T_\ast^{m+n}W) \to \mu_\ast(U)\mu_\ast(V)\mu_\ast(W).$$

Возникает вопрос:  у каких перемешивающих действий
нет эргодических гомоклинических групп?  В силу сказанного примером является перемешивающее, но не кратно перемешивающее действие Ледраппье \cite{L}.  Орициклический поток не обладает гомоклиническими элементами (кроме тождественного), что можно вывести из теоремы М. Ратнер о равномерном распределении орбит \cite{Ra}.  Интересно узнать,
бывают ли преобразования, не являющиеся элементом некоторой гомоклинической группы (в качестве кандидата можно предложить эргодический поворот окружности). 
 
\it Для перемешивающих преобразований 
 вероятностного  пространства, имеющих  положительный локальный ранг \rm (определение см. в \cite{R7}) 
можно показать, что \it гомоклиническая группа конечна, следовательно,  она не может быть эргодической. Для 
перемешивающих  преобразований ранга 1 гомоклиническая группа состоит из тождественного преобразования. \rm
Можно заметить, что любое бесконечное декартово произведение
является элементом гомоклинической группы  бернуллиевского автоморфизма с бесконечной энтропией. 
Однако  не очевидно, что существует  преобразование с нулевой  энтропией, для которого бернуллиевский автоморфизм является гомоклиническим элементом.  

М.И. Гордин поставил вопрос о существовании преобразования с нулевой энтропией и эргодической гомоклинической группой.
Положительный ответ дал Дж. Кинг в  \cite{K}.  
Наша статья содержит другие примеры: оказывается, что
все  перемешивающие  пуассоновские  надстройки, в том числе и  с нулевой энтропией\footnote{\large E. Janvresse, T. Meyerovitch, E. Roy, T. de la Rue в  работе
"Poisson suspensions and entropy for infinite transformations"
 Trans. Amer. Math. Soc. \bf 362 \rm (2010), 3069-3094,  доказали, что 
пуассоновские надстройки, индуцированные  преобразованиями 
ранга один,  имеют нулевую энтропию.}, обладают эргодической гомоклинической группой. 

Может ли преобразование с нулевой энтропией обладать эргодическим гомоклиническим   потоком?
Пример Кинга не дал ответа на этот вопрос.
Мы покажем, что найдется пуассоновская надстройка $T_\ast$ с простым 
(и, следовательно, сингулярным) спектром и  гомоклинической группой $H(T_\ast)$,  содержащей бернуллиевский поток.

Спектр пуассоновской надстройки $T_\ast$  полностью определен  спектром преобразования $T$,   поскольку  $T_\ast$  как унитарный оператор изоморфен оператору 
$exp(T)=\bigoplus_{n=0}^{\infty} T^{\odot n}$, где $T^{\odot 0}=1$ -- одномерный тождественный оператор,  $T^{\odot n}$ -- симметрическая тензорная степень оператора $T$.
Если спектр $T_\ast$ простой, то он сингулярный.
Энтропия таких надстроек равна нулю, так как 
преобразования с положительной энтропией обладают счетнократной лебеговской компонентой в спектре.

В построении подходящего бесконечного преобразования $T$
методом работы \cite{R7} добиваемся простоты спектра
$exp(T)$, при этом  обеспечиваем
наличие в группе $H(T)$ диссипативного элемента $S$.
Тогда пуассоновская надстройка $T_\ast$  будет обладать 
бернуллиевским гомоклиническим элементом $S_\ast$. 
 
В статье  будут использованы  так называемые  сидоновские конструкции 
бесконечных преобразований ранга 1.  Свойство перемешивания для них вытекает из определения,  а специфика конструкций позволяет строить диссипативные гомоклинические элементы. 
   Сидоновские  конструкции интересны и тем, что они  позволяют   получить бесконечные преобразования с простым спектром и быстрым убыванием корреляций.  Аналогичный эффект в случае вероятностного фазового пространства  впервые  был обнаружен А.А. Приходько \cite{P}.
В связи с новой теорией типичных перемешивающих преобразований ( см. недавние работы  \cite{B},\cite{T}) 
было бы интересно узнать, какие типичные свойства имеются у 
бесконечных перемешивающих преобразований.

 Автор благодарит М.И. Гордина, Э. Руа  и Ж.-П. Тувено за
вопросы, обсуждения  и полезные замечания.

% }
%\end{abstract}

\section{ Конструкции преобразований ранга один }  
Сохраняющее меру обратимое преобразование
$ T: X \to X $ пространства Лебега $ (X, \mu) $ имеет ранг 1, если  существует последовательность $ \xi_j $ измеримых разбиений  пространства 
 $ X $ таких, что 
$$ \xi_j = \{E_j, \ TE_j, \ T^2 E_j, \ \dots, T^{ h_j} E_j, \tilde {E}_j \} $$
и $ \xi_j $ стремится к разбиению на точки.
Набор "этажей" (или "блоков")  $$ E_j, TE_j, \ T^2 E_j, \ \dots, T^{ h_j} E_j $$ называется
 башней. Также башней называют дизъюнктное объединение
$ X_j  = \uu_{z = 0}^{ h_j} T^z E_j $.  В этом определении 
    мера пространства может быть  как конечной, так и бесконечной. 
Мера множества $ \tilde{E}_j = X \setminus X_j$
стремится к 0 в случае конечной меры и равна бесконечности в случае $\mu(X)=\infty$.

\it  Конструкция преобразования   ранга 1 \rm 
Фиксируем
$h_1\geq 1$ (высота башни на этапе 1), последовательность  $r_j\to\infty$ ( число колонн,  в общем случае  требуют только $r_j\geq 2$)
  и последовательность
целочисленных векторов   
$$ \bar s_j=(s_j(1), s_j(2),\dots, s_j(r_j-1),s_j(r_j)).$$  
Ниже мы даем описание конструкции, которая полностью определена параметрами $h_1$, $r_j$, $\bar s_j$.

На шаге $j=1$ задан  интервал $E_1$.
Пусть на на шаге $j$  определена  
 система (башня)  непересекающихся интервалов 
$$E_j, TE_j, T^2E_j,\dots, T^{      h_j-1}E_j,$$
причем на $E_j, TE_j, \dots, T^{      h_j-2}E_j$
пребразование $T$ является параллельным переносом.

Теперь преходим к шагу с номером $j+1$.  Представим 
интервал  $E_j$ как дизъюнктное объединение  $r_j$ интервалов  одинаковой меры, т. е.  
$$E_j=E_j^1\uu E_j^2\uu  E_j^3\uu\dots\uu E_j^{r_j}.$$
Для  $i=1,2,\dots, r_j$ рассмотрим набор  (колонну) 
$E_j^i, TE_j^i ,T^2 E_j^i,\dots, T^{h_j-1}E_j^i$
($i$-ая колонна на этапе $j$).
Над каждой колонной с номером $i$ добавим  $s_j(i)$  интервалов меры $\mu(E_j^i)$ (новые блоки в башне этапа $j+1$), получая набор  
$$E_j^i, TE_j^i, T^2 E_j^i,\dots, T^{h_j-1}E_j^i, T^{h_j}E_j^i, T^{h_j+1}E_j^i, \dots, T^{h_j+s_j(i)-1}E_j^i$$
(все множества  не пересекаются).
Обозначив $E_{j+1}= E^1_j$, для  всех $i<r_j$ положим 
$T^{h_j+s_j(i)}E_j^i = E_j^{i+1}$.
 Таким образом,  колонны складываются в одну новую башню   
$$E_{j+1}, TE_{j+1}, T^2 E_{j+1},\dots, T^{h_{j+1}-1}E_{j+1},$$
где  
 $$ 
h_{j+1} =h_jr_j +\sum_{i=1}^{r_j}s_j(i).$$

Продолжая построение, получим в качестве фазового пространства  $X$  объединение  всех интервалов и  обратимое преобразование $T$ на $X$.
Мера пространства $Х$  бесконечна, если расходится ряд   
$$\sum_j \frac {s_j(1)+s_j(2)+\dots+s_j(r_j)}{h_jr_j }.$$
 Хорошо известно, что  преобразование ранга 1  является эргодическим и имеет
простой спектр. (Простота спектра  бесконечного преобразования в общем случае не влечет за собой эргодичность.)

\section{  Сидоновские конструкции и  быстрое убывание корреляций. } 
 В этом параграфе будут рассмотрены конструкции ранга 1 со специальными  надстройками.  Такие преобразования  замечательны тем, что их свойство перемешивания легко вытекает из  определения, более того,
  перемешивание может быть быстрым. 

{\bf Сидоновские конструкции.} Пусть $r_j\to\infty$ и на каждом шаге $j$  выполнено
$$ h_j<<s_j(1)<< s_j(2)<<\dots<< s_j(r_j-1)<<s_j(r_j) \eqno (\ast).$$ 
Тогда для фиксированных $\xi_{j_0}$-измеримых множеств $A,B\subset X_{j_0}$  выполняется 
$$\mu(A\cap T^mB)\leq \mu(A)/r_j$$
для всех   $m\in [h_j, h_{j+1}],$ $j>j_0$.  
Таким образом, для всех множеств  $A,B$ конечной меры
$$\mu(A\cap T^mB)\to 0.$$

  Такая конструкция обладает \it свойством Сидона:  \it пересечение 
$X_j\cap T^mX_j$ при   $h_{j}<m\leq h_{j+1}$  может содержаться 
только в одной  колонне  башни $X_j$ \rm (будем считать сказанное  определением сидоновской конструкции).
Очевидно, что сидоновская конструкция  является перемешивающей при $r_j\to\infty$.

    \bf Oптимальные сидоновские конструкции. \rm Пусть задана последовательность  ${N_j}\to \infty$. (Отметим, что последовательность 
$N_j$ не обязана быть заданной заранее, а может определяться шаг за шагом в ходе построения преобразования.)   В  целочисленном интервале $\{1,2,\dots,N_j\}$, как известно,  найдется сидоновское множество $\bf S_j$ максимальной мощности, близкой к  $\sqrt{N_j}$. 
Напомним, что множество Сидона ${\bf S_j}$ по определению удовлетворяет уловию: при $m>0$ пересечение  ${\bf S}\cap ({\bf S}+m)$ содержит не более одного элемента.
  Обозначим
элементы  множества $\bf S_j$  через $S_j(0), S_j(1),\dots, S_j(r_j)$,
считая, что   $S_j(0)< S_j(1)<\dots <S_j(r_j)$.  
Положим $r_j=|{\bf S}_j|-1$  и
  $$s_j(i)=h_j\left(S_j(i)-S_j(i-1)-1\right), \ \ i=1, 2,\dots, r_j.$$
Такая конструкция  является сидоновской, а ее особенностью является   минимальность надстроек.  
\vspace{2mm}

\bf Теорема 3.1. \it   Для любой функции  $\psi(m)\to\infty$ (ее рост может быть сколь угодно медленным) такой, что последовательность
$\frac{\psi(m)}{\sqrt{m}}$ монотонно стремится к 0, для некоторой  оптимальной сидоновской конструкции $T$   для плотного семейства множеств $A$ (плотного
в классе всех множеств конечной меры)  выполнено условие
быстрого убывания корреляций
 $$\mu(A\cap T^mA)\leq C \frac{\psi(m)}{\sqrt{m}},$$
где константа $C$ зависит от множества $A$. \rm
\vspace{2mm}
 
Доказательство.  
На шаге  $j$  определим $r_{j}$  условием $\psi(h)\geq \sqrt{h_j}$  для всех $h>r_j^2$. Положим $N_j=r_j^2$.  Как было описано выше,
найдем $s_j(i)$, $i=1, 2,\dots, r_j,$  отвечающие оптимальной  сидоновской конструкции. Получим 
$$h_{j+1}\sim h_jN_j =  h_j r_j^2,$$  
$$\frac {\sqrt{h_j}}{\psi( h_{j+1})}\leq 1.$$
Из монотонности  $\frac{\psi(m)}{\sqrt{m}}$ для $ m\in [h_j+1, h_{j+1}]$ имеем 
$$\frac{\psi( h_{j+1})}{\sqrt{ h_{j+1}}  } \leq   \frac{\psi(m)}{\sqrt{m}}.$$

Тогда для всех $m\in [h_j+1, h_{j+1}] $ и множества $A$, состоящего из 
набора этажей некоторой башни $X_{j_0}$  выполняется 
$$\mu(A\cap T^mA)\leq \frac {\mu(A)}{r_j} \leq  C  \frac{\sqrt{h_j}}{\sqrt{ h_{j+1}}}= C
\frac{\psi( h_{j+1})}{\sqrt{ h_{j+1}}  }\frac{\sqrt{h_j}}{\psi({ h_{j+1}})}\leq C \frac{\psi(m)}{\sqrt{m}}.$$

Осталось заметить, что такие множества $A$ всюду плотны в классе всех множеств конечной меры.

\section{Кратное перемешивание, гомоклинические группы и пуассоновские надстройки. }
 Следуя  \cite{G},  определим гомоклиническую  группу
$H(T)$ преобразования   $T$,  положив 
 $$H(T)=\{ S\in Aut(X,\mu): \ T^nST^{-n}\to Id,  n\to\infty\}.$$
М. Гордин \cite{G} доказал  свойство перемешивания  автоморфизма $T$  вероятностного пространства в случае, когда группа $H(T)$ является   эргодической.   Мы усилим этот  результат.    
Автоморфизм $T$ вероятностного пространства обладает перемешиванием кратности 2, если  
для любых  измеримых множеств  $A,B,C$  при ${m,n}\to\infty$ выполнено
$$\mu (A\cap T^{m} B\cap T^{m+n} C)\to\mu (A)\mu (B)\mu(C).$$
 Аналогично определяется перемешивание кратности $n>2$.

{\bf ТЕОРЕМА 4.1.} \it  Автоморфизм вероятностного пространства  с эргодической гомоклинической группой обладает  перемешиванием любой кратности.  \rm 

Теорема очевидным образом вытекает из следующей леммы.

{\bf Лемма 4.2.} \it  Пусть $R_i$ и $T_i$ -- последовательности  автоморфизмов вероятностного  пространства.
Пусть для  автоморфизма  $S$ выполнено \rm
 $$R^{-1}_iSR_i\to Id, \ \ \ T^{-1}_iST_i\to Id.  $$

(i)\it Если для некоторой меры  
 $\nu$  на произведении $X\times X\times X$ 
для всех  измеримых множеств  $A,B,C$  выполнено 
$$\mu (A\cap R_i B\cap T_i C) \to \nu (A\times B\times C),$$ то   
$$\nu (SA\times B\times C)=\nu (A\times B\times C).$$

(ii) Пусть  $\nu (SA\times B\times C)=\nu (A\times B\times C)$ для всех $S\in H(T)$, а группа $H(T)$ эргодическая.  Если  выполнено $\nu(X\times B\times C)=\mu (B)\mu(C)$ для любых измеримых  $B,C$, то $\nu=\mu\times\mu\times\mu$.
\rm

Доказательство леммы.  (i) Имеем цепочку равенств
$$\nu(A\times B\times C)=\lim_{i\to\infty}\mu (A\cap R_i B\cap T_i C)=$$ $$%\lim_{i\to\infty}\mu (SA\cap SR_i B\cap ST_i C)
\lim_{i\to\infty} \mu (SA\cap R_iR^{-1}_iSR_i B\cap T_iT^{-1}_iST_i C)= $$ $$\lim_{i\to\infty}\mu (SA\cap R_i B\cap T_i C)=
\nu (SA\times B\times C).$$

(ii)  Мера  
$\nu$ является присоединением (joining) эргодического  действия $H(T)$  на пространстве $X_{(1)}$ и тождественного действия на пространстве $X_{(2)}\times X_{(3)}$.
Хорошо известно, что в этом случае $\nu$ является прямым произведением   проекций меры $\nu$ на указанные пространства,  т. е. $\nu=\mu \times (\mu\times\mu)$.

\bf Лемма 4.3.    \it  Если фазовое пространство имеет бесконечную меру, а преобразование $T$ перемешивает,  то группа $F$  всех преобразований $S$ с носителем конечной меры, $\mu(supp(S))<\infty$,  лежит в группе $H(T)$.  \rm

Доказательство.  Для любого множества $A$ конечной меры 
выполнено  $$\mu(A\cap supp(T^{-n}ST^n))\to 0<\infty, \ \ \ n\to\infty,$$
в силу  свойства перемешивания преобразования $T$.
Следовательно, $$T^{-n}ST^n\to Id. $$

{\bf  Пуассоновские  надстройки.}
Фиксируем пространство $(X,\mu)$ с бесконечной мерой.
На пространстве  $X_\ast$ конфигураций, т. е. счетных подмножеств пространства $X$,  вводится пуассоновская мера $\mu_\ast$. Она определена условием: для любого 
набора множеств $A_i\subset X$, $i=1,\dots, N$, конечной меры
выполнено 
$$\mu_\ast (\{x_\ast\ : \ |x_\ast\cap A_i|=n_i, \ i=1,\dots, N\})=
\prod_{i=1}^{N} \frac{e^{-\mu(A_i)}\mu(A_i)^{n_i}}
{(n_i)!}.
$$
Пусть $T$ -- автоморфизм пространства  $(X,\mu)$. Пуассоновской надстройкой
называется 
преобразование $$T_\ast (\{x_k\})=\{(Tx)_i\}.$$ 
Надстройка $T_\ast$ сохраняет вероятностную  меру
$\mu_\ast$.
Cопоставление $T\to T_\ast$ осуществляет  непрерывное вложение $Aut(X,\mu)$ в $Aut(X_\ast,\mu_\ast)$.
\vspace{3mm}

\bf Теорема 4.4. \rm (\cite{Ro}) \it Перемешивающие пуассоновские надстройки обладают перемешиванием всех порядков.\rm
\vspace{3mm}

Доказательство. Для бесконечного перемешивающего  преобразования $T$  группа $F=\{S:\ \mu(supp(S))<\infty\}$,  лежит в группе $H(T)$ (лемма 4.3).  
Пуассоновская надстройка $F_\ast$,  является эргодической группой, так как ее замыкание содержит всевозможные $T_\ast$.
Последнее есть следствие того, что замыкание $F$ содержит все автоморфизмы $T$ постранства $(X,\mu)$.  Осталось заметить, что $F_\ast$ гомоклиническая  относительно $T_\ast$ и применить теорему 4.1.

\section{ 
Пуассоновская надстройка  с сингулярным спектром и гомоклиническим бернуллиевским потоком }
Рассмотрим следующий вопрос: \it может ли автоморфизм  с нулевой энтропией обладать гомоклиническим эргодическим потоком?  \rm Мы найдем пуассоновскую надстройку с нужными свойствами.  
\vspace{2mm}

\bf  Утверждение  5.1. \it Существует пуассоновская  надстройка  с сингулярным спектром, обладающая бернуллиевским гомоклиническим потоком. \rm
\vspace{2mm}

Возьмем перемешивающий атоморфизм $R$  бесконечного пространства такой, что   $R_\ast$ имеет простой, следовательно, и  непрерывный сингулярный спектр. Примеры
имеются в работе \cite{DR}.  Для их построения  можно 
также воспользоваться сидоновскими конструкциями.
Добавляя в сидоновскую  конструкцию неперемешивающую специальную часть по аналогии с \cite{DR}, \cite{R7},  добиваемся простоты спектра оператора $exp(R)$.   Уменьшая неперемешивающую часть,  сохраняем простоту спектра  $exp(R)$, но получаем свойство перемешивания:
$$\mu(A\cap R^mB)\leq  \eps_j +\mu(A)/r_j,$$
где оценка $\eps_j\to 0$ отвечает исчезающей неперемешивающей части.
 Свойство перемешивания для таких модифицированных сидоновских $R$   не требует доказательства и вытекает из определения.

На пространстве $X=\R\times  \R^+$ рассмотрим преобразование $$T(x,y)=(x,R (y)).$$
Пусть поток $S_t$  на $X$ задан равенством
$$  S^t(x,y)=(x+\varphi(y)t, y),$$
где $\varphi >0$ и  $\varphi(y)\to 0$ при  $y\to\infty$.
Поток является диссипативным и гомоклиническим по отношению к $T$.  Первое очевидно, а второе несложно установить.  Действительно,  для функций вида  $f=\chi_{[a,b]\times [c,d]}$ выполнено
$$\|f - T^{-n}S^tT^{n} f\|_2\to 0,\ \ n\to\infty, $$
так как при фиксированном $t$  преобразование  $T^nS^tT^{-n}$  на полосе $X=\R\times  [c,d]$  близко к тождественному
преобразованию. Это видно из формулы
$$
T^{-n}S^tT^{n}(x,y)=(x+\varphi(R^n(y))t, y)
$$
и того, что в силу перемешивания для большинства 
$y\in [c,d]$  величина $R^n(y)$ велика, но 
$\varphi(y)\to 0$ при  $y\to\infty$, поэтому величина $\varphi(R^n(y))$  будет мала.
Из сказанного вытекает, что $T^{-n}S^tT^{n}\to Id$
при $n\to\infty$.

Пуассоновская надстройка $S^t_\ast$ является бернуллиевским потоком ($S^t$ -- диссипативный поток), гомоклиническим относительно $T_\ast$.
Спектр надстройки $T_\ast$ сингулярный: мера максимального спектрального типа  совпадает со спектром преобразования $R_\ast$, а  кратность спектра  надстройки $T_\ast$ бесконечна. 
Утверждение доказано.

 В следующем параграфе предлагается другой способ построения
гомоклинических действий.
\section{Пуассоновские надстройки с простым спектром и бернуллиевским гомоклиническим потоком}

Сейчас мы опишем  конструкцию $T$ ранга 1, которая имеет
в качестве одного из гомоклинических элементов диссипативное преобразование   $S$.  Напомним, что диссипативность  означает существование такого измеримого множества $Y$, что  $$X=\bigsqcup_{z\in\Z} S^zY.$$  

\bf Лемма 6.1.  \it  Пусть для некоторого автоморфизма $S$ и
 конструкции $T$ ранга 1  выполнено
$$ \min \{ \mu(S T^n T^kE_j\  | \ T^nT^kE_j)\ :
\ 0\leq k \leq h_j,\ h_j\leq n \leq h_{j+1}\} \to 1$$
при $j\to\infty$.  Тогда преобразование $S$  является гомоклиническим относительно $T$. \rm

Доказательство.  Фиксируем множество $А$, состоящее из
объединения некоторых блоков некоторой башни нашей конструкции на этапе $j_0$.  Для всех $j>j_0$  множество $A$
является объединением блоков вида $ T^kE_j$.   Пусть 
$ h_j\leq n \leq h_{j+1}$.   Так как $T^nA$ является объединением некоторых множеств вида  вида $T^nT^kE_j$,
а поледние мало отличаются от   $ST^nT^kE_j$, получаем, что $T^nA$ мало отличается от $ST^nA$.  Это означает, что 
$T^{-n}ST^nA$  асимптотически совпадает с множеством $A$.
Автоморфизм $S$ является гомоклиническим.

\bf Преобразования  $S$ и $T$, удовлетворяющие условиям леммы 6.1. \rm Пусть $T$ --  сидоновская конструкция.
При
$ h_j\leq n \leq h_{j+1}$  образ блока $ T^kE_j$ под действием 
$T^n$  будет в основном состоять из надстроек этапов $j+1$ и
$j+2$, за исключением множества, попадающего в башню
этапа $j$.  Мера этого множества не превосходит величины 
$ \frac {\mu(E_j)}{r_j}$.  

В качестве $S$  выбираем  преобразование такое, что  для некоторой последовательности  
$\eps_j\to\infty$ выполнено
$$\mu(S  T^mE_{j+1}\  | \ T^mE_{j+1})> 1-\eps_j $$ при условии, что  указанный блок $T^mE_{j+1}$  не лежит в башне  $X_j$.  
Тогда образ  блока $ T^kE_j$ под действием 
$T^n$  будет в основном состоять из таких (новых)  блоков в силу свойства Сидона  нашей конструкции.

\bf Построение диссипативного преобразования  $\bf S$. \rm
 Пространство  $X$  является дизъюнктным объединением
$X_1$  и множеств $X_j\setminus X_{j+1}$.  Последние состоят их так называемых новых блоков. Занумеруем их в порядке построения и обозначим через $D_k$.

Фиксируем последовательность $s_k\to\infty$ такую, чтобы выполнялось 
$$\sum_k \frac{\mu(D_k)}{s_k}=\infty.$$

Разобьем каждый  блок на $D_k$  на $s_k$ интервалов 
$B_k^1, B_k^2,\dots, B_k^{s_k}$ одинаковой меры.
Зададим преобразование $\tilde S$, циклически переставляющее указанные интервалы внутри блока так,  что $\tilde S^{s_k}$
является тождественным на $D_k$.

На объединении $\bigsqcup_{k}B_k^1$,  мера которого бесконечна,   рассмотрим некоторое диссипативное преобразование   $P$.
Продолжим $P$ на все пространство $X$, считая  его тождественным  вне объдинения $\bigsqcup_{k}B_k^1$.
Преобразование $S=P\tilde S$ будет диссипативным (блуждающее множество преобразования $P$ является блуждающим для $S$).
По построению, преобразование $S$  удовлетворяет свойству 
$\mu(S D_k  | \ D_k)\to 1$, $s\to\infty$. Образ  этажа $T^kE_j$ ($0\leq k\leq h_j$) под действием $T^n$ ($n>h_j$) в основном будет состоять из таких $D_{k'}$, что  значения всех таких $k'$ растут при  $j\to\infty$.
Теперь можно применить лемму 6.1.  Получим, что $S$ является
гомоклиническим относительно сидоновской конструкции  $T$. 

Отметим, что построение диссипативного  преобразования $S$ не зависит
от свойств конструкции $T$. Но для сидоновской конструкции $T$ такое преобразование оказывается гомоклиническим.  В этом случае можно применить лемму, так как  при
$ h_j\leq n \leq h_{j+1}$  образ блока $ T^kE_j$ под действием 
$T^n$  лежит  в объединении некоторых новых блоков из $X\setminus X_j$,  за исключением  асимптотически пренебрежимого множества. 

 \bf Модификация сидоновской конструкции. \rm  Рассмотрим  конструкцию $T_\eps$,
которая отличается от сидоновской тем, что на объединении  
колонн  с номерами $i, \  (1-\eps)r_j<i\leq  r_j$   не требуется выполнения  свойства Сидона.  При этом можно обеспечить
наличие слабых пределов степеней преобразования, что в свою очередь влечет за собой простоту спектра оператора
$exp(T_\eps)$. Методом, изложенным в \cite{R7},  доказывается существование  перемешивающей конструкции $T$, являющейся 
пределом таких конструкций $T_{\eps_k}$, $\eps_k\to0$, и наследующей свойство простоты спектра $exp(T)$. Аналогичная задача, но для другого класса бесконечных преобразований была решена в \cite{DR}.   Построенное  выше преобразование $S$  будет гомоклиническим относительно $T$.  Очевидно, что  $S$ вкладывается в диссипативный поток,
который является гомоклиническим относительно $T$.
Таким образом, имеет место следующее утверждение.
\vspace{4mm}

\bf Теорема  6.2. \it  Существует перемешивающая  пуассоновская надстройка с простым сингулярным спектром,  обладающая 
бернуллиевским гомоклиническим потоком. \rm
\vspace{4mm}

Замечание. Интересно узнать, всегда ли перемешивающая пуассоновская надстройка обладает бернуллиевским гомоклиническим элементом.   Вероятно, можно доказать, что   \it
существует   пуассоновская  надстройка $T_\ast$ с простым  спектром такая, что для любой
пуассоновской надстройки $R_\ast$ гомоклиническая группа
$H(T_\ast)$ содержит элемент, сопряженный с $R_\ast$ (универсальность группы $H(T_\ast)$). \rm
Для  этого  достаточно установить универсальность  $H(T)$ для модифицированной сидоновской конструкции $T$. 
\large

 E-mail: vryzh@mail.ru

\end{document}